\documentclass[10pt,a4paper,twoside]{article}
\usepackage[T1]{fontenc}
\usepackage[utf8]{inputenc}
\usepackage{lmodern}
\usepackage[textwidth=145mm,textheight=222mm,centering,headheight=14pt,headsep=7mm]{geometry}
\usepackage{amsmath,amssymb,amsthm,graphicx,booktabs,array}
\usepackage{tikz}
\usetikzlibrary{arrows.meta,calc}
\usepackage{enumitem,listings,xcolor}
\usepackage{fancyhdr,caption,cite}
\usepackage[hidelinks]{hyperref}

\hypersetup{pdftitle={Star decomposition of complete multipartite graphs},pdfauthor={Lee Juho}}
\newtheorem{theorem}{Theorem}[section]
\newtheorem{lemma}[theorem]{Lemma}
\newtheorem{proposition}[theorem]{Proposition}
\newtheorem{corollary}[theorem]{Corollary}
\theoremstyle{definition}
\newtheorem{definition}[theorem]{Definition}
\theoremstyle{remark}

\newcommand{\papertitle}[4]{\leavevmode\hbox{%
  \phantom{Decomposition of Graphs into Stars}%
  \pdfliteral page{q BT /F\pdffontname\font\space #1 Tf
    1 0 0 1 #2 #3 Tm
    [(Star) -#4 (decomposition) -#4 (of) -#4 (complete)
     -#4 (multipartite) -#4 (graphs)] TJ ET Q}%
}}
\fancypagestyle{plain}{\fancyhf{}\fancyfoot[C]{\small\thepage}}
\begin{document}
\thispagestyle{plain}
\begin{center}
\vspace*{6mm}
{\large\bfseries\papertitle{11.9552}{145.2613789}{694.672}{375}\par}
\vspace{5mm}
Juho Lee\par
\vspace{1.5mm}
{\small\itshape Korea Science Academy of KAIST, Busan, Republic of Korea\par}
{\small\texttt{24-077@ksa.hs.kr}\par}
\vspace{2mm}
{\small Research supervisor: Bumtle Kang\par}
{\small\itshape Department of Mathematics and Information Science\par}
{\small\texttt{thelokbt@ksa.kaist.ac.kr}\par}
\end{center}
\vspace{2mm}
\begin{list}{}{\leftmargin=7mm\rightmargin=7mm}
\item[]\small
We study minimum degree-sum conditions for decompositions into stars of sizes $k$ and $k+1$, where $k\ge3$. Write $\sigma(G)=\min_{uv\in E(G)}(d(u)+d(v))$, and let $c_k$ be the least integer threshold guaranteeing a $\{K_{1,k},K_{1,k+1}\}$-decomposition for every finite simple graph. Using Tarsi's sufficient condition, we prove that $K_n$ admits such a decomposition whenever $n\ge2k$. Together with the obstruction $K_{2k-1}$ of Zhao and Wu, this gives the optimal integer threshold $4k-3$ for complete graphs. For complete multipartite graphs, we establish the sufficiency of $\sigma(G)\ge9$ when $k=3$ by a finite computational verification and an extension argument. For every $k\ge3$, the complete tripartite graph $G_k=K_{k^2-k-1,k-2,1}$ satisfies $\sigma(G_k)=k^2-1$ and has no such decomposition. Consequently, $c_k\ge k^2$, even when the class is restricted to complete multipartite graphs. The optimal integer threshold in this class is therefore $9$ for $k=3$. Sufficiency of $\sigma(G)\ge k^2$ for complete multipartite graphs with $k\ge4$, and the exact values of $c_k$ for general graphs, remain open in this study.

\medskip
\noindent\textit{Keywords:} Star decomposition; degree-sum condition; complete graph; complete multipartite graph; orientation.
\end{list}
\section{Introduction}
Graph decomposition asks whether the edge set of a graph can be partitioned into subgraphs of prescribed types. We consider decompositions into stars of two consecutive sizes, $\{K_{1,k},K_{1,k+1}\}$, and investigate minimum degree-sum conditions that guarantee their existence.

Zhao and Wu~\cite{Zhao2015} proved that the condition $d(u)+d(v)\ge4k-3$ for every edge $uv$ is sufficient for the existence of a decomposition into stars of size at least $k$. Restricting the permitted sizes to exactly $k$ and $k+1$ introduces additional arithmetic constraints. Let $c_k$ denote the least integer degree-sum threshold for this restricted decomposition problem on general graphs; the precise definition is given in Section~2.

For complete graphs, we prove that $K_n$ admits the desired decomposition for every $k\ge3$ and $n\ge2k$. Zhao and Wu's obstruction $K_{2k-1}$ then shows that the optimal integer degree-sum threshold for complete graphs is $4k-3$. For complete multipartite graphs, we establish that $\sigma(G)\ge9$ is sufficient when $k=3$, using a finite verification together with a structural extension argument.

We also construct, for every $k\ge3$, the nondecomposable complete tripartite graph $G_k=K_{k^2-k-1,k-2,1}$ with $\sigma(G_k)=k^2-1$. Thus a threshold guaranteeing decomposition is at least $k^2$, even within the complete multipartite class, and in particular $c_k\ge k^2$. For $k\ge4$, this rules out the sufficiency of $4k-3$. For $k=3$, it proves that the optimal integer threshold for complete multipartite graphs is $9$. Sufficiency of $k^2$ for general $k$ requires further investigation.

\section{Background and Previous Work}
All graphs in this paper are finite and simple. We follow the terminology of Bondy and Murty~\cite{Bondy2008}. For a graph $G=(V(G),E(G))$, write $d(v)$ for the degree of $v$ and $N(v)$ for its neighborhood. For $S\subseteq V(G)$, $G[S]$ denotes the subgraph induced by $S$: its vertices are $S$, and its edges are exactly those edges of $G$ with both endpoints in $S$.

\begin{definition}\label{def:graph_types}
The following three graph types will be used throughout.
\begin{enumerate}[leftmargin=*,label=(\roman*),itemsep=3pt,topsep=3pt]
\item \textbf{Complete graph.} The graph $K_n$ has $n$ vertices and an edge between every pair of distinct vertices.
\item \textbf{Complete multipartite graph.} The graph $K_{n_1,\ldots,n_r}$ has a partition $V_1,\ldots,V_r$ of its vertex set into independent sets of sizes $n_1,\ldots,n_r$, with every possible edge between distinct parts.
\item \textbf{Star.} A star of size $m$, denoted by $K_{1,m}$ or $S_m$, consists of a center joined to $m$ leaves.
\end{enumerate}
\end{definition}
\begin{figure}[htbp]
\centering
\tikzset{every node/.style={circle,draw,minimum size=4mm,fill=black!8}}
\begin{tikzpicture}[scale=0.8]
                \foreach \i in {1,...,4} {
                    \node (\i) at ({\i*90+45}:1.2) {};
                }
                \draw (1) -- (2); \draw (1) -- (3); \draw (1) -- (4);
                \draw (2) -- (3); \draw (2) -- (4);
                \draw (3) -- (4);
            \end{tikzpicture}\hspace{10mm}
\begin{tikzpicture}[scale=0.8]
                
                \node (a1) at (-0.8, 0.7) {};
                \node (a2) at (-0.8, -0.7) {};
                
                \node (b1) at (0.8, 1) {};
                \node (b2) at (0.8, 0) {};
                \node (b3) at (0.8, -1) {};
                
                \foreach \i in {1,2} {
                    \foreach \j in {1,2,3} {
                        \draw (a\i) -- (b\j);
                    }
                }
            \end{tikzpicture}\hspace{10mm}
\begin{tikzpicture}[scale=0.8]
                \node (center) at (0,0) {};
                \foreach \i in {1,...,4} {
                    \node (\i) at ({(\i-1)*90}:1.2) {};
                    \draw (center) -- (\i);
                }
            \end{tikzpicture}
\caption{Examples in the order of Definition~\ref{def:graph_types}: $K_4$, $K_{2,3}$, and $K_{1,4}$.}
\end{figure}
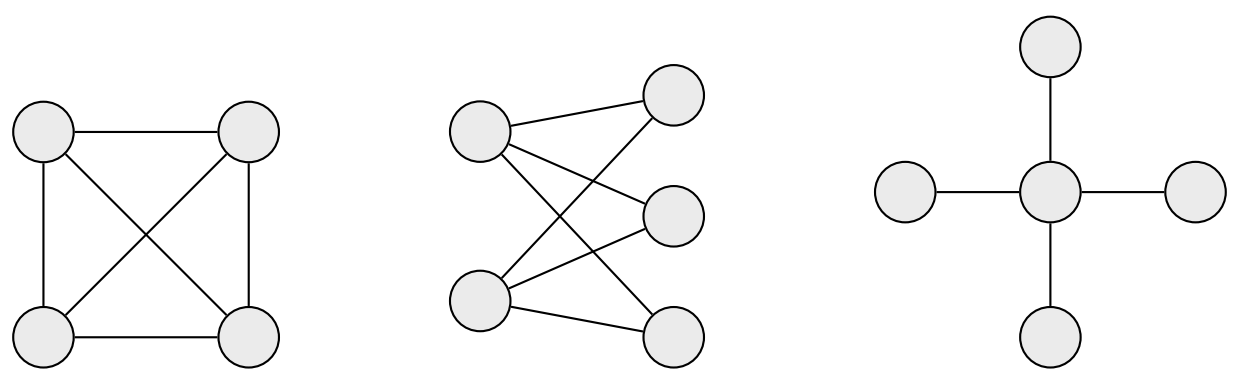
\begin{definition}
A decomposition of $G$ is a family of pairwise edge-disjoint subgraphs whose edge sets have union $E(G)$. It is an $\mathcal F$-decomposition if each member is isomorphic to a graph in $\mathcal F$.
\end{definition}

\begin{definition}
For a graph with at least one edge, define its minimum degree sum by
\[
\sigma(G)=\min_{uv\in E(G)}\bigl(d(u)+d(v)\bigr).
\]
Thus $\sigma(G)\ge c$ is equivalent to $d(u)+d(v)\ge c$ for every edge $uv$. An edgeless graph has the empty decomposition and is treated separately whenever $\sigma$ is used.
\end{definition}
For a fixed integer $k\ge3$, let $c_k$ be the least integer $c$ such that every finite simple graph with at least one edge and $\sigma(G)\ge c$ admits a $\{K_{1,k},K_{1,k+1}\}$-decomposition. We discuss thresholds for complete multipartite graphs directly, without introducing a separate symbol.

\subsection{Known sufficient conditions}
\begin{theorem}[Zhao and Wu~\cite{Zhao2015}]\label{thm:zhao}
Let $G$ be a graph and $k\ge3$ an integer. If $d(u)+d(v)\ge4k-3$ for every edge $uv$, then $G$ admits a star decomposition in which every star has size at least $k$.
\end{theorem}
Every integer $q\ge k(k-1)$ can be expressed as $q=ak+b(k+1)$ with $a,b\in\mathbb Z_{\ge0}$. Applying Theorem~\ref{thm:zhao} with parameter $k(k-1)$ and splitting the resulting stars gives the following consequence.
\begin{corollary}\label{corollary:k_k+1}
If $G$ has at least one edge and $\sigma(G)\ge4k(k-1)-3$, then $G$ admits a $\{K_{1,k},K_{1,k+1}\}$-decomposition.
\end{corollary}
Together with the lower bound proved in Corollary~\ref{cor:k2_lower_bound}, this yields
\[
k^2\le c_k\le4k(k-1)-3.
\]
The inequality $c_k\ge k^2$ rules out every smaller integer threshold. It does not show that $c_k=k^2$: proving sufficiency of $\sigma(G)\ge k^2$ requires a separate argument.

\subsection{Orientations and complete graphs}
Our approach replaces star decomposition by a constraint on outdegrees. Classical results of Landau~\cite{Landau1953} and Tarsi~\cite{Tarsi1981} provide background for this viewpoint.
\begin{theorem}[Landau~\cite{Landau1953}]
A nondecreasing sequence of nonnegative integers $(s_1,\ldots,s_n)$ is the outdegree sequence of an orientation of $K_n$ if and only if
\[
\sum_{i=1}^{j}s_i\ge\binom{j}{2}\quad(1\le j<n),
\qquad \sum_{i=1}^{n}s_i=\binom{n}{2}.
\]
\end{theorem}
Tarsi~\cite{Tarsi1979} studied decompositions of complete multigraphs into stars of equal size. His later result~\cite{Tarsi1981} gives the following sufficient condition for stars of prescribed sizes.
\begin{theorem}[Tarsi~\cite{Tarsi1981}]\label{thm:tarsi}
Let $m_1,\ldots,m_t$ be positive integers. The complete graph $K_n$ admits a decomposition into stars of these sizes if
\[
\sum_{i=1}^{t}m_i=\binom n2,
\qquad 1\le m_i\le n/2\quad(1\le i\le t).
\]
\end{theorem}
\begin{lemma}[Zhao and Wu~\cite{Zhao2015}, Lemma 2.4]\label{lem:zhao_complete_obstruction}
For every integer $k\ge3$, $K_{2k-1}$ has no decomposition into stars each containing at least $k$ edges.
\end{lemma}
In particular, $K_{2k-1}$ has no $\{K_{1,k},K_{1,k+1}\}$-decomposition. Combining this obstruction with our sufficient condition for $n\ge2k$ determines the optimal integer threshold within the complete graph class.

The results below address three objectives: decomposability of complete graphs for all $n\ge2k$; sufficiency of $\sigma(G)\ge9$ for complete multipartite graphs when $k=3$; and a family of complete tripartite obstructions establishing the lower bound $k^2$ for every $k\ge3$.

\section{An Equivalent Orientation Condition and Computational Methods}
\subsection{An equivalent outdegree condition}
\begin{lemma}[Orientation criterion]\label{lem:orientation}
A finite simple graph $G$ admits a $\{K_{1,k},K_{1,k+1}\}$-decomposition if and only if it has an orientation such that, at every vertex $v$,
\begin{equation}\label{eq:orientation}
d^+(v)=ka(v)+(k+1)b(v),\qquad a(v),b(v)\in\mathbb Z_{\ge0}.
\end{equation}
\end{lemma}
\begin{proof}
Given a decomposition, orient each star from its center toward its leaves. Each edge belongs to exactly one star, so it receives exactly one direction. If $a(v)$ and $b(v)$ count the stars of sizes $k$ and $k+1$ centered at $v$, respectively, then~\eqref{eq:orientation} follows.

Conversely, suppose such an orientation exists. At each vertex $v$, partition its outgoing edges into $a(v)$ groups of size $k$ and $b(v)$ groups of size $k+1$. Each group forms a star centered at $v$. Every edge has exactly one tail and is therefore used exactly once. These stars form the required decomposition.
\end{proof}

For $k=3$, the allowed outdegrees are $\{0,3,4,6,7,8,\ldots\}$; equivalently, the forbidden outdegrees are $1$, $2$, and $5$.

\subsection{Mixed-integer formulation and exploratory experiments}
We used SageMath mixed-integer linear programming to generate graphs and test decomposability. The decision constraints implement Lemma~\ref{lem:orientation}. For each edge $uv$, binary variables satisfy $x_{uv}+x_{vu}=1$, and the vertex constraints are
\[
\sum_{u\in N(v)}x_{vu}=ka(v)+(k+1)b(v),
\qquad a(v),b(v)\in\mathbb Z_{\ge0}.
\]
In our exploratory experiment, the generator used orders $9\le n\le24$, the condition $\sigma(G)\ge9$, and the additional constraint $|E(G)|\ge n$. Its objective had independently chosen edge coefficients from $[-1,1]$. The generator is given in Appendix~\ref{app:generation}, and the general orientation-based decision test is given in Appendix~\ref{app:general}. The procedure does not sample uniformly from all qualifying graphs; in particular, the edge-count constraint excludes every forest.

In our exploratory experiment, $500{,}000$ trials were performed and no counterexample was found. The experiment provides evidence within the generated sample, rather than an exhaustive verification or a proof of universal sufficiency. The exploratory computations used SageMath.

Appendix~\ref{app:random_driver} presents the random-experiment loop, with the trial count set to $500{,}000$. It selects the graph order, displays progress, counts successful decisions and counterexamples, and stores counterexample graphs in a list. Failed graph generation and unresolved solver decisions are counted separately. Appendix~\ref{app:multipartite_driver} presents the recursive partition enumeration and complete multipartite test loop, configured with $k=3$ and orders $6$ through $10$. This loop tests every partition and classifies confirmed counterexamples by whether their minimum degree sum is below $4k-3$.

The mixed-integer decision tests in Appendix~\ref{app:general} and Appendix~\ref{app:structured} use GLPK (GNU Linear Programming Kit), a solver for linear and mixed-integer programming. They return \texttt{True} after a successful solve, \texttt{False} only when GLPK explicitly reports infeasibility, and \texttt{None} for other solver exceptions. An unresolved instance must not be counted as nondecomposable. Edgeless graphs return \texttt{True}.

The finite cases in Section~\ref{sec:multipartite3} and the counterexample tables in Section~\ref{sec:computational} were additionally checked by exhaustive integer orientation search using Python~3.12.14, without an external optimization solver. This check used the orientation criterion in Lemma~\ref{lem:orientation} independently of the MILP tests in Appendices~\ref{app:general} and~\ref{app:structured}. The $120$ small complete multipartite cases are summarized in Table~\ref{tab:small}; Appendix~\ref{app:multipartite_driver} gives their partition-enumeration loop and MILP tests.

\section{Decomposition Results and Lower Bounds}
\subsection{Complete graphs}
Since every vertex of $K_n$ has degree $n-1$, we have $\sigma(K_n)=2(n-1)$. For integer $n\ge2$,
\[
\sigma(K_n)\ge4k-3\quad\Longleftrightarrow\quad n\ge2k.
\]
\begin{theorem}[Complete graph decomposition]\label{thm:complete}
Let $k\ge3$ and $n\ge2$ be integers. If $\sigma(K_n)\ge4k-3$, equivalently $n\ge2k$, then $K_n$ admits a $\{K_{1,k},K_{1,k+1}\}$-decomposition.
\end{theorem}
\begin{proof}
By Theorem~\ref{thm:tarsi}, it suffices to choose star sizes $m_i\in\{k,k+1\}$ with sum $\binom n2$ and $m_i\le n/2$.

First suppose $n\ge2k+2$. Then $k+1\le n/2$, so either permitted star size satisfies the bound in Theorem~\ref{thm:tarsi}. Moreover,
\[
\binom n2\ge(2k+1)(k+1)=2k^2+3k+1>k^2-k,
\]
since the difference between the last two expressions is $k^2+4k+1>0$. The Frobenius coin theorem for the coprime integers $k$ and $k+1$ therefore gives nonnegative integers $a,b$ such that $\binom n2=ak+b(k+1)$. Taking $a$ stars of size $k$ and $b$ stars of size $k+1$ satisfies all conditions of Theorem~\ref{thm:tarsi}.

For $n=2k$ or $n=2k+1$, the larger size $k+1$ exceeds $n/2$, so we use only stars of size $k$. The identities
\[
\binom{2k}{2}=k(2k-1),\qquad
\binom{2k+1}{2}=k(2k+1)
\]
show that the edge count is divisible by $k$ in both cases, and $k\le n/2$. Theorem~\ref{thm:tarsi} yields a $\{K_{1,k}\}$-decomposition, which is also a $\{K_{1,k},K_{1,k+1}\}$-decomposition.
\end{proof}

\begin{corollary}[Optimal integer threshold for complete graphs]\label{cor:complete_optimal_threshold}
For every $k\ge3$, the least integer $c$ such that every complete graph $K_n$, $n\ge2$, with $\sigma(K_n)\ge c$ admits a $\{K_{1,k},K_{1,k+1}\}$-decomposition is $4k-3$.
\end{corollary}
\begin{proof}
Theorem~\ref{thm:complete} proves sufficiency. By Lemma~\ref{lem:zhao_complete_obstruction}, $K_{2k-1}$ has no decomposition even into stars all of size at least $k$, whereas
\[
\sigma(K_{2k-1})=2(2k-2)=4k-4.
\]
Thus no integer $c\le4k-4$ guarantees the desired decomposition for all complete graphs.
\end{proof}
Because $\sigma(K_n)$ is even, the conditions $\sigma(K_n)\ge4k-3$ and $\sigma(K_n)\ge4k-2$ are equivalent. Optimality here refers to the \emph{least integer threshold}. This result is restricted to complete graphs and does not determine $c_k$ for general graphs.

\subsection{Complete multipartite graphs when \texorpdfstring{$k=3$}{k=3}}\label{sec:multipartite3}
\begin{lemma}\label{lem:extension}
Let $G$ be a complete multipartite graph with parts $V_1,\ldots,V_r$, and let $H=G[\bigcup_{i\in I}V_i]$ for an index set $I$. If $H$ has a $\{K_{1,3},K_{1,4}\}$-decomposition and $|V(H)|\ge6$, then $G$ has such a decomposition.
\end{lemma}
\begin{proof}
Let $W_0=V(H)$ and enumerate $V(G)\setminus W_0$ as
$u_1,\ldots,u_t$. For $j=0,\ldots,t$, put
\[
W_j=W_0\cup\{u_1,\ldots,u_j\},\qquad H_j=G[W_j].
\]
Thus $H_0=H$ and $H_t=G$. We use mathematical induction on $j$ to prove that every $H_j$ admits a $\{K_{1,3},K_{1,4}\}$-decomposition.

\smallskip
\noindent\textbf{Base case ($j=0$).}
The graph $H_0=H$ admits such a decomposition by hypothesis. This also proves the assertion immediately if $t=0$.

\smallskip
\noindent\textbf{Inductive step.}

Suppose $1\le j\le t$ and $H_{j-1}$ has such a decomposition. Define
\[
F_j=\{u_jv\in E(G):v\in W_{j-1}\},\qquad q_j=|F_j|.
\]
Since $W_0$ is a union of entire partite sets and $u_j\notin W_0$,
the vertex $u_j$ belongs to none of the selected parts. It is therefore
adjacent to every vertex of $W_0$, and $q_j\ge|W_0|\ge6$.
The integers $3$ and $4$ are coprime and their Frobenius number is
$3\cdot4-3-4=5$. Hence $q_j=3a_j+4b_j$ for some
$a_j,b_j\in\mathbb Z_{\ge0}$.
Partition $F_j$ into $a_j$ sets of three edges and $b_j$ sets of four
edges. Each set forms a star centered at $u_j$, since $G$ is simple.
Moreover,
\[
E(H_j)=E(H_{j-1})\mathbin{\dot\cup}F_j.
\]
Adding these stars to the decomposition of $H_{j-1}$ gives a
decomposition of $H_j$, completing the inductive step.

\smallskip
By mathematical induction, every $H_j$, $0\le j\le t$, admits the required decomposition. In particular, $H_t=G$ does.
\end{proof}
\begin{lemma}\label{lem:subset_sum}
Let $n_1,\ldots,n_r$ be positive integers, each at most $5$. If $\sum_i n_i\ge11$, some index set $I$ satisfies $6\le\sum_{i\in I}n_i\le10$.
\end{lemma}
\begin{proof}
For $j=0,\ldots,r$, define
\[
s_0=0,\qquad s_j=\sum_{i=1}^{j}n_i,
\]
and let
\[
j_0=\min\{j\in\{1,\ldots,r\}:s_j\ge6\}.
\]
This index exists because $s_r\ge11$. By minimality,
$s_{j_0-1}<6$, so $s_{j_0-1}\le5$, since the partial sums are integers.
Since $n_{j_0}\le5$, it follows that
\[
6\le s_{j_0}=s_{j_0-1}+n_{j_0}\le5+5=10.
\]
Taking $I=\{1,\ldots,j_0\}$ gives
$6\le\sum_{i\in I}n_i\le10$.
The choice is made by indices, so distinct parts of equal size remain
distinct.
\end{proof}
\begin{theorem}[Complete multipartite decomposition for $k=3$]\label{thm:multipartite3}
Every complete multipartite graph $G$ with at least one edge and $\sigma(G)\ge9$ admits a $\{K_{1,3},K_{1,4}\}$-decomposition. Edgeless graphs also admit such a decomposition.
\end{theorem}
\begin{proof}
The edgeless case has the empty decomposition. Otherwise, write $G=K_{n_1,\ldots,n_r}$, where $n_1\ge n_2\ge\cdots\ge n_r\ge1$ and $r\ge2$. Let $V_i$ be the part of size $n_i$ for each $i\in\{1,\ldots,r\}$, and put $N=\sum_{i=1}^r n_i$. A vertex in $V_i$ has degree $N-n_i$, so
\[
\sigma(G)=2N-(n_1+n_2)\ge9.
\]
If $n_1\ge6$, the edgeless graph $G[V_1]$ has at least $6$ vertices and Lemma~\ref{lem:extension} applies. We may therefore assume $n_1\le5$.

The degree condition implies $N\ge(9+n_1+n_2)/2\ge11/2$, hence $N\ge6$. If $N\le10$, select all parts. If $N\ge11$, use Lemma~\ref{lem:subset_sum} to select entire parts whose total size is between $6$ and $10$. Let $H$ be their induced graph. If $H$ is decomposable, Lemma~\ref{lem:extension} finishes the proof. The degree-sum condition need not pass from $G$ to $H$, so the exceptional small graphs require separate treatment.

There are $120$ complete multipartite graphs of orders $6$ through $10$, up to isomorphism, corresponding to integer partitions of those orders. Of these, $94$ have all part sizes at most $5$. The exhaustive integer orientation search found exactly four nondecomposable graphs:
\[
K_{5,1},\quad K_{4,1,1},\quad K_{5,2},\quad K_{5,1,1}.
\]
All four lie in the restricted set of $94$ graphs. Thus $116$ of the full set and $90$ of the restricted set are decomposable. The counts by order are as follows.
\begin{table}[htbp]
\centering
\begin{tabular}{crrl}
\toprule
Order & All partitions & Largest part $\le5$ & Exceptions\\
\midrule
6 & 11 & 10 & $K_{5,1},K_{4,1,1}$\\
7 & 15 & 13 & $K_{5,2},K_{5,1,1}$\\
8 & 22 & 18 & None\\
9 & 30 & 23 & None\\
10 & 42 & 30 & None\\
\midrule
Total & 120 & 94 & $K_{5,1},K_{4,1,1},K_{5,2},K_{5,1,1}$\\
\bottomrule
\end{tabular}
\caption{Finite verification for complete multipartite graphs when $k=3$.}\label{tab:small}
\end{table}
These finite cases were also checked by exhaustive integer orientation search using Python~3.12.14, without an external optimization solver. The general and structured SageMath decision routines appear in Appendices~\ref{app:general} and~\ref{app:structured}. The finite enumeration loop in Appendix~\ref{app:multipartite_driver} uses the general decision routine. This check uses the orientation criterion in Lemma~\ref{lem:orientation} independently of the MILP tests.
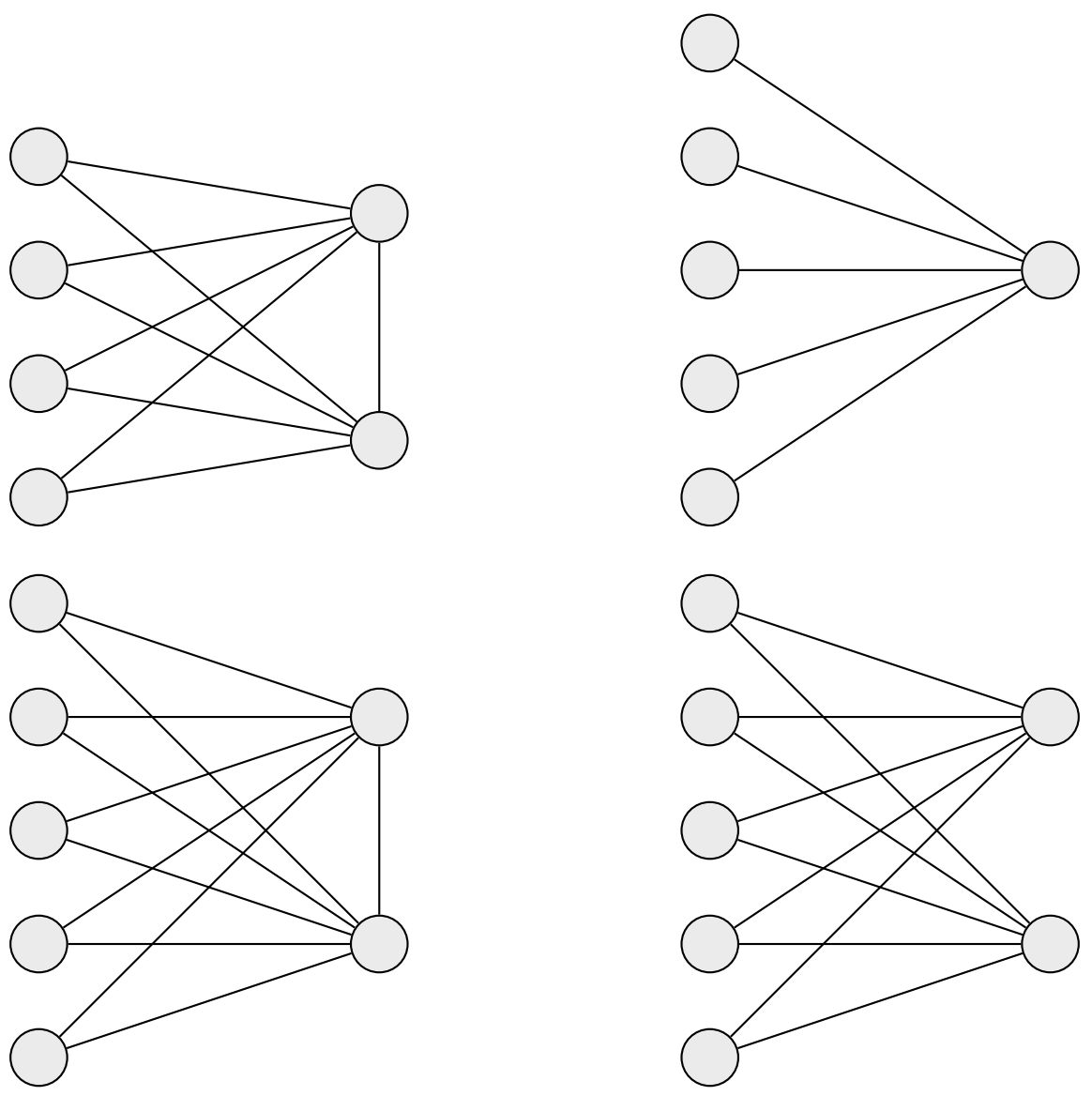
\begin{figure}[htbp]
\centering
\tikzset{every node/.style={circle,draw,minimum size=4mm,fill=black!8}}
\begin{tikzpicture}[scale=0.8]
                
                \foreach \i in {1,...,4} {
                    \node (a\i) at (-1.5, {2.5-\i}) {};
                }
                
                \node (b1) at (1.5, 1) {};
                
                \node (c1) at (1.5, -1) {};

                \foreach \i in {1,...,4} {
                    \draw (a\i) -- (b1);
                    \draw (a\i) -- (c1);
                }
                \draw (b1) -- (c1); 
            \end{tikzpicture}\hspace{18mm}
\begin{tikzpicture}[scale=0.8]
                
                \foreach \i in {1,...,5} {
                    \node (a\i) at (-1.5, {3-\i}) {};
                }
                
                \node (b1) at (1.5, 0) {};

                \foreach \i in {1,...,5} {
                    \draw (a\i) -- (b1);
                }
            \end{tikzpicture}\par\vspace{3mm}
\begin{tikzpicture}[scale=0.8]
                
                \foreach \i in {1,...,5} {
                    \node (a\i) at (-1.5, {3-\i}) {};
                }
                
                \node (b1) at (1.5, 1) {};
                
                \node (c1) at (1.5, -1) {};

                \foreach \i in {1,...,5} {
                    \draw (a\i) -- (b1);
                    \draw (a\i) -- (c1);
                }
                \draw (b1) -- (c1); 
            \end{tikzpicture}\hspace{18mm}
\begin{tikzpicture}[scale=0.8]
                
                \foreach \i in {1,...,5} {
                    \node (a\i) at (-1.5, {3-\i}) {};
                }
                
                \node (b1) at (1.5, 1) {};
                \node (b2) at (1.5, -1) {};

                \foreach \i in {1,...,5} {
                    \draw (a\i) -- (b1);
                    \draw (a\i) -- (b2);
                }
            \end{tikzpicture}
\caption{The four small nondecomposable graphs: $K_{4,1,1}$, $K_{5,1}$, $K_{5,1,1}$, and $K_{5,2}$.}
\end{figure}
The four graphs $K_{4,1,1}$, $K_{5,1}$, $K_{5,1,1}$, and $K_{5,2}$ have minimum degree sums $7$, $6$, $8$, and $7$, respectively. Hence $G$ cannot equal an exceptional $H$. At least one unselected part remains; let its size be $m$, where $1\le m\le5$.
\begin{enumerate}[leftmargin=*,label=(\roman*)]
\item If $H=K_{4,1,1}$ and $1\le m\le4$, the enlarged graph $K_{4,1,1,m}$ has order $7$ through $10$ and is not exceptional, so it is decomposable. If $m=5$, it contains the decomposable graph $K_{5,4}$ as a union of entire parts, and Lemma~\ref{lem:extension} applies.
\item If $H=K_{5,1}$ and $2\le m\le4$, the graph $K_{5,1,m}$ has order $8$ through $10$ and is not exceptional. If $m=5$, apply Lemma~\ref{lem:extension} to $K_{5,5}$. If $m=1$, the enlarged graph is $K_{5,1,1}$, whose minimum degree sum is still only $8$. Another part must remain, and the next case applies.
\item If $H=K_{5,1,1}$ and $1\le m\le3$, the enlarged graph has order $8$ through $10$ and is not exceptional. For $m=4,5$, apply Lemma~\ref{lem:extension} to $K_{5,4}$ or $K_{5,5}$, respectively.
\item If $H=K_{5,2}$ and $1\le m\le3$, the enlarged graph has order $8$ through $10$ and is not exceptional. For $m=4,5$, again apply Lemma~\ref{lem:extension} to $K_{5,4}$ or $K_{5,5}$.
\end{enumerate}
In every case we obtain a decomposable induced graph on at least $6$ vertices consisting of entire parts. Lemma~\ref{lem:extension} then extends the decomposition to $G$.
\end{proof}

\subsection{Explicit computational counterexamples}\label{sec:computational}
The computational counterexamples below motivated the general construction in Section~\ref{sec:general_obstruction}. The structured-graph decision functions are given in Appendix~\ref{app:structured}, and the partition-enumeration and degree-sum classification procedure is given in Appendix~\ref{app:multipartite_driver}.
Among the examples for $k=3,4,5$ are $K_{5,1,1}$, $K_{11,2,1}$, and $K_{19,3,1}$, with minimum degree sums $8,15,24$, respectively. These graphs motivated the family in Proposition~\ref{prop:multipartite_counterexample}. The corresponding thresholds are at least $9,16,25$, so $c_3\ge9$, $c_4\ge16$, and $c_5\ge25$.

For $k=4$, the previously considered condition was $\sigma(G)\ge13=4k-3$. Six nondecomposable complete multipartite graphs meeting this condition were found; see Table~\ref{tab:counterexamples_k4}. For example, $K_{10,1,1}$ has $12$ vertices and minimum degree sum $13$, but no orientation with every outdegree of the form $4a+5b$.
\begin{table}[htbp]
\centering
\begin{tabular}{c|c|c|c}
\hline
\textbf{Order} & \textbf{Parts} & $\boldsymbol{\sigma(G)}$ & \textbf{Decomposable} \\ \hline
12 & $K_{10, 1, 1}$ & 13 & No \\
13 & $K_{11, 2}$ & 13 & No \\
13 & $K_{11, 1, 1}$ & 14 & No \\
13 & $K_{10, 1, 1, 1}$ & 15 & No \\
14 & $K_{11, 3}$ & 14 & No \\
14 & $K_{11, 2, 1}$ & 15 & No \\ \hline
\end{tabular}
\caption{Nondecomposable graphs for $k=4$ satisfying $\sigma(G)\ge13$.}
\label{tab:counterexamples_k4}
\end{table}
For $k=5$, the analogous condition was $\sigma(G)\ge17$. The search found the $27$ nondecomposable examples listed in Table~\ref{tab:counterexamples_k5}.
\begin{table}[htbp]
\centering
\begin{tabular}{c|c|c|c}
\hline
\textbf{Order} & \textbf{Parts} & $\boldsymbol{\sigma(G)}$ & \textbf{Decomposable} \\ \hline
15 & $K_{12, 1, 1, 1}$ & 17 & No \\
15 & $K_{11, 2, 1, 1}$ & 17 & No \\
15 & $K_{11, 1, 1, 1, 1}$ & 18 & No \\
16 & $K_{14, 1, 1}$ & 17 & No \\
16 & $K_{13, 2, 1}$ & 17 & No \\
16 & $K_{13, 1, 1, 1}$ & 18 & No \\
16 & $K_{12, 2, 2}$ & 18 & No \\
16 & $K_{12, 2, 1, 1}$ & 18 & No \\
16 & $K_{12, 1, 1, 1, 1}$ & 19 & No \\
17 & $K_{14, 3}$ & 17 & No \\
17 & $K_{14, 2, 1}$ & 18 & No \\
17 & $K_{13, 4}$ & 17 & No \\
17 & $K_{13, 3, 1}$ & 18 & No \\
17 & $K_{13, 2, 2}$ & 19 & No \\
17 & $K_{13, 2, 1, 1}$ & 19 & No \\
17 & $K_{13, 1, 1, 1, 1}$ & 20 & No \\
18 & $K_{14, 4}$ & 18 & No \\
18 & $K_{14, 3, 1}$ & 19 & No \\
20 & $K_{19, 1}$ & 20 & No \\
20 & $K_{18, 1, 1}$ & 21 & No \\
21 & $K_{19, 2}$ & 21 & No \\
21 & $K_{19, 1, 1}$ & 22 & No \\
21 & $K_{18, 1, 1, 1}$ & 23 & No \\
22 & $K_{19, 3}$ & 22 & No \\
22 & $K_{19, 2, 1}$ & 23 & No \\
23 & $K_{19, 4}$ & 23 & No \\
23 & $K_{19, 3, 1}$ & 24 & No \\ \hline
\end{tabular}
\caption{Nondecomposable graphs for $k=5$ satisfying $\sigma(G)\ge17$.}
\label{tab:counterexamples_k5}
\end{table}

\clearpage
\subsection{A general obstruction and the quadratic lower bound}\label{sec:general_obstruction}
The computational examples suggest a family valid for every $k\ge3$. Here $\sigma(G)=k^2-1$ specifies the \emph{minimum} over edges; it does not assert equal degree sums on all edges.
\begin{proposition}\label{prop:multipartite_counterexample}
For every integer $k\ge3$, the complete tripartite graph
\[
G_k=K_{k^2-k-1,\,k-2,\,1}
\]
satisfies $\sigma(G_k)=k^2-1$ and has no $\{K_{1,k},K_{1,k+1}\}$-decomposition.
\end{proposition}
\begin{proof}
Put $L=k^2-k-1$, and denote the parts by $A,B,C=\{z\}$, with sizes $L,k-2,1$. Their vertex degrees are
\[
d(a)=k-1\ (a\in A),\quad d(b)=k^2-k\ (b\in B),\quad d(z)=k^2-3.
\]
The possible degree sums are
\[
\begin{array}{c|c}
\text{Edge type}&\text{Degree sum}\\\hline
A\text{--}B&k^2-1\\
A\text{--}C&k^2+k-4\\
B\text{--}C&2k^2-k-3
\end{array}
\]
For $k\ge3$, the last two values are at least $k^2-1$, and edges between $A$ and $B$ exist. Thus $\sigma(G_k)=k^2-1$.

First, $L$ cannot be written as $\alpha k+\beta(k+1)$ for nonnegative integers $\alpha,\beta$. Such an expression would imply $\beta\equiv k-1\pmod k$, so $\beta\ge k-1$, whereas
\[
\beta(k+1)\ge(k-1)(k+1)=k^2-1>L,
\]
a contradiction.

Suppose a decomposition exists. By Lemma~\ref{lem:orientation}, orient every star from its center to its leaves. Each edge belongs to exactly one star, so this defines an orientation of all of $G_k$. If $\alpha_v$ and $\beta_v$ count stars of sizes $k$ and $k+1$ centered at $v$, then
\[
d^+(v)=\alpha_vk+\beta_v(k+1),\qquad
\alpha_v,\beta_v\in\mathbb Z_{\ge0}.
\]
Thus no vertex can have outdegree $L$. Each vertex of $A$ has degree $k-1<k$, so its outdegree must be $0$. Every edge between $B\cup C$ and $A$ is consequently directed toward $A$.

Each $b\in B$ already has $L$ outgoing edges toward $A$ and only one other incident edge, namely $bz$. Its outdegree is therefore $L$ or $L+1$. Since $L$ is forbidden, every edge $bz$ must be oriented $b\to z$. The only outgoing edges of $z$ are then its $L$ edges toward $A$, so $d^+(z)=L$, again impossible. No decomposition exists. The forced directions are summarized in Figure~\ref{fig:forced}.
\end{proof}

\begin{figure}[htbp]
\centering
\begin{tikzpicture}[>=Stealth, thick, every node/.style={font=\small}]
\node[draw,rounded corners,align=center,minimum width=31mm,minimum height=16mm] (A) at (0,0) {$A$: $L$ vertices\\$d^+(a)=0$};
\node[draw,rounded corners,align=center,minimum width=31mm,minimum height=16mm] (B) at (7,1.5) {$B$: $k-2$ vertices\\$d^+(b)=L+1$};
\node[draw,circle,align=center,minimum size=14mm] (Z) at (7,-1.5) {$z$};
\draw[->,blue!65!black] (B.west) -- node[above,sloped] {(1) all $B$--$A$ edges} (A.north east);
\draw[->,blue!65!black] (Z.west) -- node[below,sloped] {(1) all $z$--$A$ edges} (A.south east);
\draw[->,red!70!black] (B.south) -- node[right,align=left] {(2) each $b\to z$\\since $L$ is forbidden} (Z.north);
\node[align=center] at (3.5,-3) {(3) $d^+(z)=L$ is forbidden: contradiction.};
\end{tikzpicture}
\caption{Forced directions in Proposition~\ref{prop:multipartite_counterexample}, where $L=k^2-k-1$. Each arrow to $A$ represents all edges between the indicated parts, rather than a single edge.}
\label{fig:forced}
\end{figure}
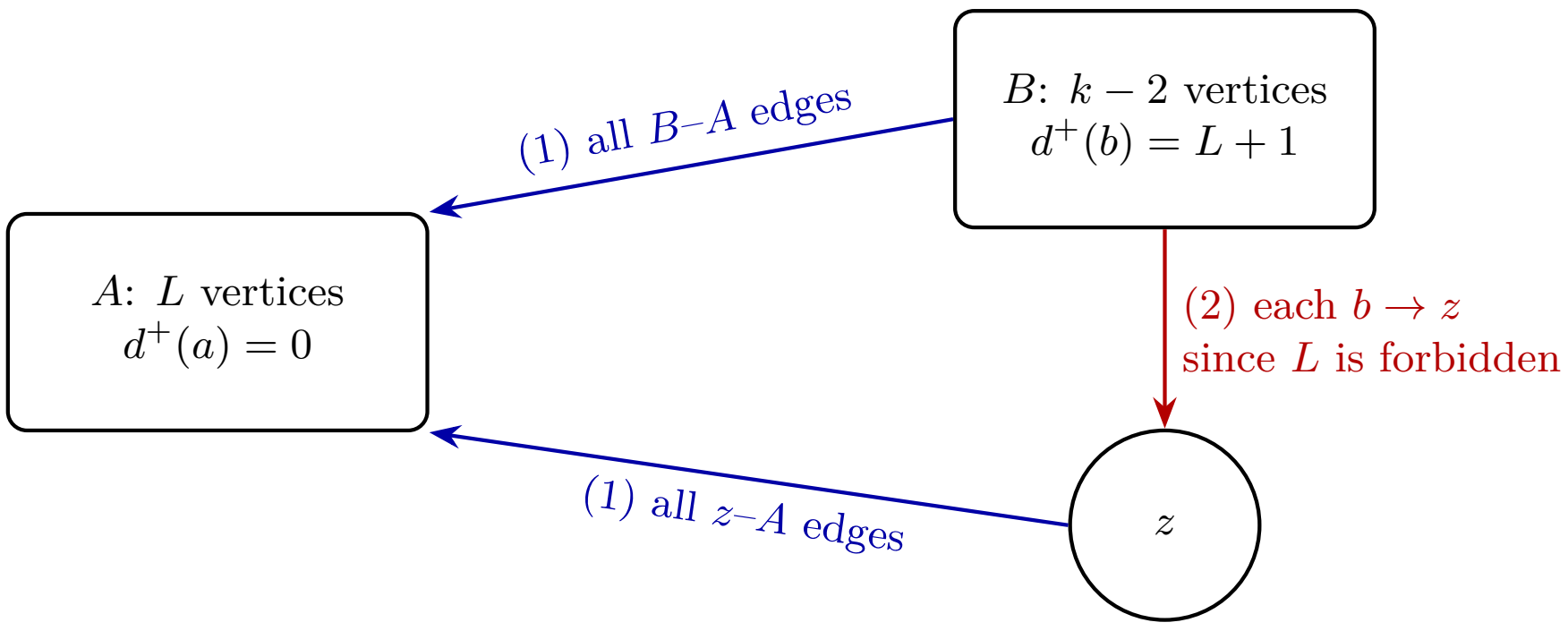
\begin{corollary}\label{cor:k2_lower_bound}
For every $k\ge3$, any integer threshold guaranteeing a $\{K_{1,k},K_{1,k+1}\}$-decomposition for all complete multipartite graphs is at least $k^2$. Consequently, $c_k\ge k^2$. For $k=3$, the optimal integer threshold for complete multipartite graphs is $9$.
\end{corollary}
\begin{proof}
The obstruction in Proposition~\ref{prop:multipartite_counterexample} has minimum degree sum $k^2-1$, so every sufficient integer threshold must exceed this value. Complete multipartite graphs are a subclass of general graphs, giving $c_k\ge k^2$. For $k=3$, combine the obstruction with Theorem~\ref{thm:multipartite3}.
\end{proof}

The obstruction family proves a lower bound for every $k\ge3$. Since $k^2-(4k-3)=(k-1)(k-3)>0$ for $k\ge4$, the older candidate $4k-3$ is insufficient even for complete multipartite graphs. At $k=3$, the lower bound is attained in this class. Whether $\sigma(G)\ge k^2$ is sufficient for $k\ge4$, or for general graphs, remains unresolved here.

\section{Discussion}
Throughout this discussion, decomposition means a $\{K_{1,k},K_{1,k+1}\}$-decomposition.

For complete graphs, the arithmetic representation of the edge count and Tarsi's sufficient condition establish decomposability whenever $k\ge3$ and $n\ge2k$. Combining this with Zhao and Wu's obstruction $K_{2k-1}$ yields the exact least integer degree-sum threshold $4k-3$. The sufficient condition comes from the present argument, while the matching obstruction is an existing result. Since degree sums of complete graphs are even, thresholds $4k-3$ and $4k-2$ impose the same condition; the former is the least integer threshold. This class-specific result does not determine $c_k$ for general graphs.

For complete multipartite graphs and $k=3$, the sufficient condition relies on a finite verification for small graphs and an extension argument using entire parts. Together with $\sigma(K_{5,1,1})=8$ and nondecomposability, it gives the exact threshold $9$. No counterexample was found in our $500{,}000$ exploratory trials on generated general graphs. This observation motivates further study but does not prove $c_3=9$.

The general lower bound arises from the tripartite structure of $G_k=K_{k^2-k-1,k-2,1}$. Vertices of degree $k-1$ cannot be star centers. Their incident edges force orientations on the other vertices, eventually producing an outdegree that is not a nonnegative combination of $k$ and $k+1$. This yields $\sigma(G_k)=k^2-1$ and hence the lower bound $k^2$, even in the complete multipartite class.

The central remaining question is whether $\sigma(G)\ge k^2$ guarantees decomposition for complete multipartite graphs. The current results prove necessity of this threshold but do not establish sufficiency. Orientations exploiting multipartite structure, degree-control techniques, and network-flow formulations may help settle this question and improve bounds for general graphs.

\section{Conclusions and Further Questions}
Starting from Zhao and Wu's star-decomposition results, we studied minimum degree-sum conditions for decompositions into $K_{1,k}$ and $K_{1,k+1}$. For complete graphs, Tarsi's theorem gives a decomposition whenever $k\ge3$ and $\sigma(K_n)\ge4k-3$. The obstruction $K_{2k-1}$ shows that $4k-3$ is the optimal integer threshold in this class. For complete multipartite graphs, $\sigma(G)\ge9$ is sufficient when $k=3$.

For every $k\ge3$, the complete tripartite graph $K_{k^2-k-1,k-2,1}$ has minimum degree sum $k^2-1$ and no desired decomposition. Thus every sufficient threshold for the complete multipartite class is at least $k^2$, and for general graphs
\[
k^2\le c_k\le4k(k-1)-3.
\]
The complete multipartite threshold is exactly $9$ when $k=3$. For $k\ge4$, $4k-3$ is insufficient even within this class. Our exploratory random experiment provides exploratory evidence and must be distinguished from a proof of sufficiency.

Further work concerns the following questions.
\begin{enumerate}[leftmargin=*,label=(\roman*)]
\item Does $\sigma(G)\ge k^2$ guarantee a $\{K_{1,k},K_{1,k+1}\}$-decomposition for every complete multipartite graph when $k\ge4$?
\item Does $\sigma(G)\ge9$ suffice for general graphs when $k=3$? More generally, can the upper bound on $c_k$ be reduced toward $k^2$?
\item Can network-flow and degree-control methods characterize orientations whose outdegrees lie in $\{ak+b(k+1):a,b\in\mathbb Z_{\ge0}\}$?
\end{enumerate}

\clearpage
\appendix
\section{Python (SageMath) Code}
Appendices~\ref{app:general}--\ref{app:structured} give the decision and graph-generation routines; Appendices~\ref{app:random_driver} and~\ref{app:multipartite_driver} add execution drivers for the random experiment and finite multipartite searches. The decision routines use GLPK as specified in the source. The integer $n$ in \texttt{can\_decompose\_Kn\_Knplus1} is the star-size parameter denoted by $k$ in the text.
\subsection{General decomposition test}\label{app:general}
\begin{lstlisting}
from sage.numerical.mip import MixedIntegerLinearProgram, MIPSolverException

def can_decompose_Kn_Knplus1(G, n):
    if G.size() == 0:
        return True
    p = MixedIntegerLinearProgram(maximization=False, solver="GLPK")
    
    x = p.new_variable(binary=True)
    
    a = p.new_variable(integer=True, nonnegative=True)
    b = p.new_variable(integer=True, nonnegative=True)
    
    for u, v in G.edges(labels=False):
        p.add_constraint(x[u, v] + x[v, u] == 1)
        
    for v in G.vertices():
        out_degree_v = sum(x[v, u] for u in G.neighbors(v))
        p.add_constraint(out_degree_v == n * a[v] + (n+1) * b[v])
        
    try:
        p.solve()
        return True
    except MIPSolverException as exc:
        # False only for confirmed GLPK infeasibility; otherwise unknown.
        message = str(exc).lower()
        if message in ("glpk: problem has no feasible solution",
                       "glpk: problem has no integer feasible solution"):
            return False
        return None
\end{lstlisting}
\subsection{Graph generation}\label{app:generation}
\begin{lstlisting}
from sage.all import Graph
from sage.numerical.mip import MixedIntegerLinearProgram, MIPSolverException
import random

def create_custom_graph_ilp(n, k):
    p = MixedIntegerLinearProgram(maximization=False)
    x = p.new_variable(binary=True)
    
    d = {i: sum(x[i, j] for j in range(n) if i != j) for i in range(n)}
    
    for i in range(n):
        for j in range(i + 1, n):
            p.add_constraint(x[i, j] == x[j, i])
            
    for i in range(n):
        for j in range(i + 1, n):
            p.add_constraint(d[i] + d[j] >= k * x[i, j])
            
    p.add_constraint(sum(x[i, j] for i in range(n) for j in range(i + 1, n)) >= n)
    
    p.set_objective(sum(random.uniform(-1, 1) * x[i, j] for i in range(n) for j in range(i + 1, n)))
    
    try:
        p.solve()
        G = Graph(n)
        for i in range(n):
            for j in range(i + 1, n):
                if p.get_values(x[i, j]) > 0.5:
                    G.add_edge(i, j)
        return G
    except MIPSolverException:
        return None
\end{lstlisting}
\subsection{Tests for complete and complete multipartite graphs}\label{app:structured}
\begin{lstlisting}
from sage.all import *
from sage.numerical.mip import MIPSolverException

def can_decompose_star_multipartite(parts, k_star):
    num_groups = len(parts)
    if num_groups < 2:
        return True

    total_vertices = sum(parts)
    sum_of_squares = sum(p**2 for p in parts)
    total_edges = (total_vertices**2 - sum_of_squares) // 2
    
    if total_edges < k_star**2 - k_star:
        representable = False
        for a_sum in range((total_edges // k_star) + 1):
            if (total_edges - a_sum * k_star) % (k_star + 1) == 0:
                representable = True
                break
        if not representable:
            return False
        
    vertices = []
    for group_id, size in enumerate(parts):
        for _ in range(size):
            vertices.append(group_id)
            
    n = len(vertices)
    p = MixedIntegerLinearProgram(maximization=False, solver="GLPK")
    
    d = p.new_variable(binary=True)
    a = p.new_variable(integer=True, nonnegative=True)
    b = p.new_variable(integer=True, nonnegative=True)
    
    for u in range(n):
        for v in range(u + 1, n):
            if vertices[u] != vertices[v]:
                p.add_constraint(d[u, v] + d[v, u] == 1)
                
    for u in range(n):
        max_deg = total_vertices - parts[vertices[u]]
        
        p.set_max(a[u], max_deg // k_star)
        p.set_max(b[u], max_deg // (k_star + 1))
        
        out_degree_u = sum(d[u, v] for v in range(n) if vertices[u] != vertices[v])
        
        p.add_constraint(out_degree_u == k_star * a[u] + (k_star + 1) * b[u])
        
    try:
        p.solve(log=0) 
        return True
    except MIPSolverException as exc:
        # False only for confirmed GLPK infeasibility; otherwise unknown.
        message = str(exc).lower()
        if message in ("glpk: problem has no feasible solution",
                       "glpk: problem has no integer feasible solution"):
            return False
        return None
        

def is_representable(target, n):
    if target >= n**2 - n:
        return True
    for a in range(target // n + 1):
        if (target - a * n) % (n + 1) == 0:
            return True
    return False


def solve_complete_graph(N, n):
    if N % 2 == 1:
        if is_representable((N - 1) // 2, n): return True
    else:
        if is_representable(N // 2, n) and is_representable(N // 2 - 1, n): return True

    valid_scores = [x for x in range(N) if is_representable(x, n)]
    if not valid_scores:
        return False
    
    memo = {}
    def dfs(k, current_sum, last_score):
        if k == N:
            return current_sum == N * (N - 1) // 2
        
        if current_sum < k * (k - 1) // 2:
            return False
            
        state = (k, current_sum, last_score)
        if state in memo:
            return memo[state]
            
        for score in valid_scores:
            if score >= last_score:
                if dfs(k + 1, current_sum + score, score):
                    memo[state] = True
                    return True
                    
        memo[state] = False
        return False

    return dfs(0, 0, 0)
\end{lstlisting}

\subsection{Execution driver for the random experiment}\label{app:random_driver}
\begin{lstlisting}
# Random experiment loop.
import random
from generate_random_graph import create_custom_graph_ilp
from decomposition import can_decompose_Kn_Knplus1

n, k = 15, 9
trial = 500000
progress_unit = 1
t, f = 0, 0
u, generation_failed = 0, 0
li = []

possible_n = list(range(9, 25))
def progress(i,total, unit):
    return int(i*100/total/unit)
    
for i in range(trial):
    n = random.choice(possible_n)
    if progress(i,trial,progress_unit) != progress(i-1,trial,progress_unit): print(f"{progress(i,trial,progress_unit)*progress_unit}% complete, T:{t}, F:{f}")
    G = create_custom_graph_ilp(n, k)
    if G is None:
        generation_failed += 1
        continue
    result = can_decompose_Kn_Knplus1(G, 3)
    if result is True: t += 1
    elif result is None:
        u += 1
    else:
        f += 1
        li.append(G)
print(f"T: {t}, F: {f}, Unknown: {u}, Generation failed: {generation_failed}")
\end{lstlisting}
\subsection{Execution driver for complete multipartite searches}\label{app:multipartite_driver}
\begin{lstlisting}
# Complete multipartite graph enumeration.
from sage.all import graphs
from decomposition import can_decompose_Kn_Knplus1


def min_adjacent_degree_sum(G):
    # The minimum is undefined when the graph has no edges.
    if G.size() == 0:
        return None
    return min(G.degree(u) + G.degree(v)
               for u, v in G.edges(labels=False))

def get_partitions_with_repetition(n, start=1):
    if n == 0:
        return [[]]
    
    results = []
    for i in range(start, n + 1):
        sub_partitions = get_partitions_with_repetition(n - i, i)
        
        for p in sub_partitions:
            results.append([i] + p)
            
    return results

p = 3

false_graph = dict()
for n in range(6,11):
    
    vertex_li = get_partitions_with_repetition(n)
    
    t,f,df = 0, 0, 0
    u = 0
    li = []
    
    unit = 5
    total = len(vertex_li)
    def progress(i,total, unit):
        return int(i*100/total/unit)
    graphs_to_save = dict()
    graphs_to_save['df'] = []
    graphs_to_save['f'] = []
    i = 0
    print(f'Analysis started: {total} graphs')
    for v_li in vertex_li:
        G = graphs.CompleteMultipartiteGraph(v_li)
        result = can_decompose_Kn_Knplus1(G,p)
        if result is True: t += 1
        elif result is None:
            u += 1
        else:
            if min_adjacent_degree_sum(G) < 4*p-3:
                df += 1
                graphs_to_save['df'].append(G)
            else:
                f += 1
                graphs_to_save['f'].append(G)
            li.append(G)
        # if progress(i,total,unit) != progress(i-1,total,unit): print(f"{progress(i,total,unit)*unit}% complete, T: {t}, F: {f}, DF: {df}")
        i += 1
    false_graph[n] = graphs_to_save
    print(f"Order: {n}, T: {t}, F: {f}, DF: {df}, Unknown: {u}")
\end{lstlisting}

\begin{thebibliography}{9}\small
\bibitem{Bondy2008} J. A. Bondy and U. S. R. Murty, \textit{Graph Theory} (Springer, 2008).
\bibitem{Landau1953} H. G. Landau, On dominance relations and the structure of animal societies: III. The condition for a score structure, \textit{Bull. Math. Biophys.} \textbf{15}(2) (1953) 143--148. \url{https://doi.org/10.1007/BF02476378}.
\bibitem{Tarsi1979} M. Tarsi, Decomposition of complete multigraphs into stars, \textit{Discrete Math.} \textbf{26}(3) (1979) 273--278. \url{https://doi.org/10.1016/0012-365X(79)90034-7}.
\bibitem{Tarsi1981} M. Tarsi, On the decomposition of a graph into stars, \textit{Discrete Math.} \textbf{36}(3) (1981) 299--304. \url{https://doi.org/10.1016/S0012-365X(81)80025-8}.
\bibitem{Zhao2015} Y. Zhao and B. Wu, Star decomposition of graphs, \textit{Discrete Math. Algorithms Appl.} \textbf{7}(2) (2015) 1550016. \url{https://doi.org/10.1142/S1793830915500160}.
\end{thebibliography}
\end{document}